# Modularity Classes and Boundary Effects in Multivariate Stochastic Dominance

Luciano A. Perez*[1]
Departamento de Matemática, Escuela de Ingeniería y Gestión, Instituto Tecnológico de Buenos Aires
Av. Eduardo Madero 399
C1106ACD – Ciudad de Buenos Aires – Argentina
*luperez@itba.edu.ar

**Abstract:** Hadar & Russell (1974) and Levy & Paroush (1974) presented sufficient conditions for multivariate stochastic dominance when the distributions involved are continuous with compact support. Further generalizations involved either independence assumptions (Sacarsini (1988)) or the introduction of new concepts like "correlation increasing transformation" (Epstein & Tanny (1980), Tchen (1980), Mayer (2013)). In this paper, we present a direct proof that extends the original results to the general case where the involved distributions are only assumed to have compact support. This result has in turn proven useful for statistical tests of dominance without the assumption of absolute continuity.

The first section introduces several concepts used throughout the paper. In the second section we recall the classic result as presented in Atkinson & Bourguignon (1982), with a slightly lighter proof using the general integration by parts formula for Lebesgue-Stieljes integrals in $\mathbb{R}^n$. In the third section we present our proof of the general result, using Riemman-Stieljes partial sums in a direct fashion that helps to clarify the role of modularity conditions and boundary effects in the sufficiency of the conditions. The last section discusses the relevance of the result and concludes.

---

[1] This paper is part of a PHD research program, funded with CONICET scholarships and under the direction of Dr. Juana Z. Brufman.

# 1. Introduction

The concept of stochastic dominance was introduced in 1969 by Hadar & Russell[2] in the context of risk analysis. In essence, stochastic dominance is a mathematical rule for ordering univariate or multivariate stochastic prospects, according to their effect on the expected value of some objective functional.

Combined with the use of von Neumann-Morgenstern utility functions over lotteries, the concept of stochastic dominance has been useful in the areas of portfolio analysis and risk management for over 40 years[3].

Notwithstanding its applications in investment theory, one of the main reasons for the popularity of the notion of stochastic dominance are its applications to the study of economic inequality, beginning with the seminal paper by Atkinson in 1970[4].

Both in finance and economic inequality the basic concept of stochastic dominance in the univariate case goes as follows: for a given measurable function $u(x)$ we can define a functional $W_u$ over the space $\Im$ of cumulative distribution functions[5] with support in the interval [0,1]:

$$W_u(F) = \int_0^1 u(t)dF(t) \qquad (1.1)$$

Let $F, G \in \Im$. If for some class $\mathcal{U}$ of functions we have $W_u(F) \geq W_u(G)$ for every member $u$ in the class we say that $F$ stochastically dominates $G$ over the class $\mathcal{U}$.

The most relevant class $\mathcal{J}$ is that of increasing $\mathcal{C}^1$ functions, i.e. those with $u'(x) \geq 0$. If $F$ stochastically dominates $G$ over this class, we say that there is *first order stochastic dominance* of $F$ over $G$ and we denote this by $F\ SD_1\ G$.

Within $\mathcal{J}$ the most relevant subclass is that of concave functions, naturally associated with risk aversion. If $F$ dominates $G$ over the class of $\mathcal{C}^2$ functions with $u'(x) \geq 0$ and $u''(x) \leq 0$ we say that there is second order stochastic dominance of $F$ over $G$ and we denote this by $F\ SD_2\ G$.

---

[2] Hadar, J. & Russell, W. R. (1969). Rules for Ordering Uncertain Prospects. *Amer. Econ. Rev. 59*, 25-34.

[3] Levy, H. (2006). *Stochastic Dominance. Investment Decision Making under Uncertainty.* Second Edition, New York: Springer.

[4] Atkinson, A. (1970). On the measurement of Inequality. *Journal of Economic Theory, 2, 244-263.*

[5] This can be immediately generalized to a space of *càdlàg* functions with arbitrary (common) compact support.

Necessary conditions for univariate first order stochastic dominance $F\ SD_1\ G$ are trivially obtained by taking $u(t) = 1_{(-\infty, x]}(t)$ which yields $F(x) \leq G(x)$ for $x \in [0,1]$.

The fact that this condition is also *sufficient* for first order dominance is not trivial at all and was first proved by Hadar & Russell in their 1969 paper. Perhaps the easiest proof is the one using quantile couplings, by showing that first order dominance implies point wise dominance of certain couplings and then using the fact that the functions $u \in \mathcal{I}$ are increasing[6].

This is in fact a special case of Strassen's Theorem for partially ordered polish spaces, as shown by Perez[7].

Second order univariate stochastic dominance was proven by Hadar & Russell in their 1969 paper to be equivalent to:

$$\int_0^z F(x)dx \leq \int_0^z G(x)dx \qquad (1.2)$$

for every $z \in [0,1]$. This naturally leads to a general definition of arbitrary order stochastic dominance. We start constructing a sequence of operators over the space $\mathfrak{J}$ of cumulative distribution functions:

$$\mathcal{S}_1(z,F) = F(z) \quad \mathcal{S}_{j+1}(z,F) = \int_0^z \mathcal{S}_j(t,F)dt \qquad (1.3)$$

We say that there is $j$ order stochastic dominance of $F$ over $G$ (noted as $F\ SD_j\ G$) if:

$$\mathcal{S}_j(z,F) \leq \mathcal{S}_j(z,G) \qquad (1.4)$$

for all $z \in [0,1]$.

It is clear that the nested construction of the integral operators implies a sequential hierarchy of dominance, that is: $F\ SD_j\ G$ implies $F\ SD_{j+k}\ G$ for every $k$.

This formulation of stochastic dominance arises as a natural extension of Hadar & Russell's results and its importance lies on the fact that it translates the condition of dominance from an infinite set of inequalities (one for each function $u$) into a condition stated with a single inequality involving cumulative distribution functions $F$ and $G$.

This in turn allows an easier statistical testing of stochastic dominance, using versions of the operators obtained by plug-in of the empirical distribution

---

[6] Thorisson, H. (2000). *Coupling, Stationarity, and Regeneration. Probability and its Applications,* New York: Springer.

[7] Perez, L. (2015). Tests de Dominancia Estocástica en base a Estadísticos de Kolmogorov-Smirnov Multivariados, con Aplicaciones al Estudio de la Desigualdad Económica Multidimensional. PhD. Dissertation, UBA.

functions into the integral operators, as done by McFadden[8] and by Davidson & Duclos[9].

## 2. Hadar & Russell (1974) and Levy & Paroush (1974) results

Given two random vectors *(X$_1$,Y$_1$)* and *(X$_2$,Y$_2$)* with respective cdfs *F$_1$(s,t)* and *F$_2$(s,t)* an analogue condition for first order stochastic dominance *F$_1$ SD$_1$ F$_2$* would be *F$_1$(s,t)≤F$_2$(s,t)*. But for the multivariate case this condition, although necessary, is not sufficient.

Hadar & Russell[10] and Levy & Paroush[11] found sufficient conditions for bivariate first and second order stochastic dominance, involving integral operators as in the univariate case. The key point here is to restrict our attention to dominance over certain classes of functions with properties that make integration by parts easier in the bivariate case.

As we are using bivariate Lebesgue-Stieljes integration by parts, we are going to need conditions on second order derivatives even for first order dominance conditions. The classification of bivariate functions based on the sign of second order derivatives gives rise to the definition of *modularity classes.*

**Definition 2.1.** (Modularity Classes). Let $f:\mathbb{R}^n \to \mathbb{R}$ be a differentiable function. We say that $f$ is *supermodular* if $\frac{\partial f}{\partial x_i}(x)$ is non decreasing in $x_j$ for every pair of indexes $i \neq j$. If $f \in \mathcal{C}^2$ then $f$ is supermodular if and only if $\frac{\partial^2 f}{\partial x_i \partial x_j}(x) \geq 0$. If $\frac{\partial^2 f}{\partial x_i \partial x_j}(x) \leq 0$ then $f$ is *submodular.* The class of supermodular functions is noted by $\mathcal{M}^+$ and that of submodular functions by $\mathcal{M}^-$. □

Modularity classes may seem restrictive but most of the usual classes of utility functions, e.g. Cobb-Douglas functions, have definite modularity. In the context of microeconomics, supemodularity of a bivariate utility function would be associated with strategic complementarity of input goods[12].

---

[8] McFadden, D. (1989). Testing for Stochastic Dominance. In Fomby, Th. B & Tae Kun Seo (Eds.), *Studies in the Economics of Uncertainty. In Honor of Josef Hadar*, New York: Springer.

[9] Davidson, R. & Duclos, J-Y. (2000). Statistical Inference for Stochastic Dominance and for the Measurement of Poverty and Inequality. *Econometrica, 68 (6)*, 1435-1464.

[10] Hadar, J. & Russell, W. R. (1974). Stochastic Dominance in Choice under Uncertainty. In M. S. Balch, D. L. McFadden & Y. Wu (Eds.), *Essays on Economic Behavior under Uncertainty* (pp.133-150). North Holland, Amsterdam.

[11] Levy, H. & Paroush, J. (1974). Toward Multivariate Efficiency Criteria. *Journal of Economic Theory, 7*, 129-142.

Generalizing the univariate concept we'll say that a random $\mathbb{R}^n$ element $Z$ stochastically dominates another random $\mathbb{R}^n$ element $W$ if for all $\phi:\mathbb{R}^n \to \mathbb{R}$ in a certain class of Borel functions we have $E(\phi(Z)) \geq E(\phi(W))$.

The following theorem gathers both Hadar & Russell's (1974) and Levy & Paroush's (1974) results on sufficient conditions for absolutely continuous bivariate distributions. Its proof follows the idea of that of Atkinson & Bourguignon[13] but using the general integration by parts formula in $\mathbb{R}^n$ which allows for a lighter, more elegant proof.

**Theorem 2.2.** (First Order Bivariate Stochastic Dominance, Absolutely Continuous Case). Let $(X_1, Y_1)$ and $(X_2, Y_2)$ be random vectors with cdfs $F_1$ and $F_2$ respectively. Let's suppose in addition that both distributions are absolutely continuous with respect to Lebesgue measure and have common support with compact closure $\bar{U} \subseteq \mathbb{R}^2$. Then:

1) If $F_1(s,t) \leq F_2(s,t)$ for all $(s,t) \in U$, then $E(\phi(X_1, Y_1)) \geq E(\phi(X_2, Y_2))$ holds for every $\mathcal{C}^2$ function $\phi \in \mathcal{M}^-$ which is increasing in each argument.
2) Let's denote the marginal distribution of $X_1$ in the first vector by $F_1^X$ and by $F_1^Y$ the marginal cdf for $Y_1$. $F_2^X$ and $F_2^Y$ are analogous for the second vector. We define $K_1(s,t) = F_1^X(s) + F_1^Y(t) - F_1(s,t)$ and its analogue $K_2(s,t)$ for every $(s,t) \in U$.
If we assume that $F_1^X(s) \leq F_2^X(s)$, $F_1^Y(t) \leq F_2^Y(t)$ and $K_1(s,t) \leq K_2(s,t)$ $\forall (s,t) \in U$, then it holds $E(\phi(X_1, Y_1)) \geq E(\phi(X_2, Y_2))$ holds for every $\mathcal{C}^2$ function $\phi \in \mathcal{M}^+$ which is increasing in each argument.

**Proof.** As distributions are absolutely continuous we have density functions $\frac{\partial^2 F_1}{\partial x \partial y}(s,t) = f_1(s,t)$ and $\frac{\partial^2 F_2}{\partial x \partial y}(s,t) = f_2(s,t)$ and then:

$$E(\phi(X_1, Y_1)) = \int_U \phi(s,t) f_1(s,t) ds dt \qquad (2.1)$$

Let's recall the integration by parts formula for a pair of functions $u, v \in \mathcal{C}^1(\bar{U})$ where $U \subseteq \mathbb{R}^n$ is and open bounded set (without loss of generality we can assume $\bar{U} = [0,1] \times [0,1]$):

$$\int_U \frac{\partial u}{\partial x_i} v \, dx = -\int_U u \frac{\partial v}{\partial x_i} dx + \int_{\partial U} uv \breve{n}^i dS \qquad (2.2)$$

---

[12] Topkis, D. (1998). *Supermodularity and Complementarity*. Princeton University Press.

[13] Atkinson, A. B. & Bourguignon, F. (1982). The Comparison of Multi-Dimensioned Distributions of Economic Status, *The Review of Economic Studies, 49 (2)*, April, 183-201.

where $\partial U$ is the boundary of $U$ oriented with external normal $\breve{n}$ and $\breve{n}^i = \breve{n} \cdot e_i$ is its component in direction $x_i$.

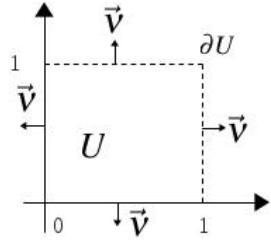

**Fig.2.1**. U and its oriented boundary.

Using this formula on our expectation we have:

$$E\big(\phi(X_1, Y_1)\big) = \int_U \phi(s,t) \frac{\partial^2 F_1}{\partial x \partial y}(s,t) ds dt = -\int_U \frac{\partial \phi}{\partial x}(s,t) \frac{\partial F_1}{\partial y}(s,t) ds dt + \int_{\partial U} \phi(s,t) \frac{\partial F_1}{\partial y}(s,t) \breve{n}^x dS = I_1 + I_2 \qquad (2.3)$$

Applying integration by parts one more time we have:

$$I_1 = \int_U \frac{\partial^2 \phi}{\partial x \partial y}(s,t) F_1(s,t) ds dt - \int_{\partial U} \frac{\partial \phi}{\partial x}(s,t) F_1(s,t) \breve{n}^y ds dt = I_3 + I_4 \qquad (2.4)$$

Taking into account the orientation of $\breve{n}$ along $\partial U$ and the fact that $F_1(0,t) = F_1(s,0) = 0$ we have:

$$I_2 = \int_{\partial U} \phi(s,t) \frac{\partial F_1}{\partial y}(s,t) \breve{n}^x dS = -\int_0^1 \frac{\partial \phi}{\partial y}(1,t) F_1^Y(t) dt + \phi(1,1) \qquad (2.5)$$

Using the same tools we have:

$$I_4 = -\int_0^1 \frac{\partial \phi}{\partial x}(s,1) F_1^X(s) ds \qquad (2.6)$$

Gathering all calculations we have:

$$E\big(\phi(X_1, Y_1)\big) = \int_U \frac{\partial^2 \phi}{\partial x \partial y}(s,t) F_1(s,t) ds dt - \int_0^1 \frac{\partial \phi}{\partial y}(1,t) F_1^Y(t) dt - \int_0^1 \frac{\partial \phi}{\partial x}(s,1) F_1^X(s) ds + \phi(1,1) = A_1 + B_1 + C_1 + \phi(1,1) \qquad (2.7)$$

Obvious analogues for $F_2$ yield:

$$E\big(\phi(X_2, Y_2)\big) = A_2 + B_2 + C_2 + \phi(1,1) \qquad (2.8)$$

Now if $\phi \in \mathcal{M}^-$ then $\frac{\partial^2 \phi}{\partial x \partial y}(s,t) \leq 0$, and as $F_1(s,t) \leq F_2(s,t)$ it follows that

$$\frac{\partial^2 \phi}{\partial x \partial y}(s,t) F_2(s,t) \leq \frac{\partial^2 \phi}{\partial x \partial y}(s,t) F_1(s,t) \qquad (2.9)$$

and then $A_2 \leq A_1$. On the other hand, taking limits on the inequality $F_1(s,t) \leq F_2(s,t)$ gives us $F_1^X(s) \leq F_2^X(s)$ and $F_1^Y(t) \leq F_2^Y(t)$. The fact that $\frac{\partial \phi}{\partial x} \geq 0$ and $\frac{\partial \phi}{\partial y} \geq 0$ then accounts for $B_2 \leq B_1$ and $C_2 \leq C_1$. This proves part 1 of Theorem 2.2.

To prove part 2 we first use the identity $\frac{\partial \phi}{\partial y}(1,t) = \int_0^1 \frac{\partial^2 \phi}{\partial x \partial y}(s,t)ds + \frac{\partial \phi}{\partial y}(0,t)$ to reformulate $B_1$ as:

$$B_1 = -\int_U \frac{\partial^2 \phi}{\partial x \partial y}(s,t)F_1^Y(t)dsdt - \int_0^1 \frac{\partial \phi}{\partial y}(0,t)F_1^Y(t)dt \tag{2.10}$$

Analogously:

$$C_1 = -\int_U \frac{\partial^2 \phi}{\partial x \partial y}(s,t)F_1^X(t)dsdt - \int_0^1 \frac{\partial \phi}{\partial x}(s,0)F_1^X(s)ds \tag{2.11}$$

which, replacing in (2.7) and using the definition of $K_1(s,t)$ leads to:

$$E(\phi(X_1,Y_1)) = -\int_U \frac{\partial^2 \phi}{\partial x \partial y}(s,t)K_1(s,t)dsdt - \int_0^1 \frac{\partial \phi}{\partial x}(s,0)F_1^X(s)ds - \int_0^1 \frac{\partial \phi}{\partial y}(0,t)F_1^Y(t)dt \tag{2.12}$$

Using $\frac{\partial \phi}{\partial x} \geq 0$ and $\frac{\partial \phi}{\partial y} \geq 0$ and the same inequalities for cdfs as before, yields $E(\phi(X_1,Y_1)) \geq E(\phi(X_2,Y_2))$ as wanted. □

For second order dominance we need to define higher order modularity classes.

**Definition 2.3.** (Higher Order Modularity Classes). Let $U \subseteq \mathbb{R}^n$ and Let $f \in \mathcal{C}^2(U)$. We say that $f \in \mathcal{M}^{--}$ if $f \in \mathcal{M}^-$ and:

$$\begin{cases} \frac{\partial^2 f}{\partial x^2}, \frac{\partial^2 f}{\partial y^2} \leq 0 \\ \frac{\partial^3 f}{\partial x^2 \partial y}, \frac{\partial^3 f}{\partial x \partial y^2} \geq 0 \\ \frac{\partial^4 f}{\partial x^2 \partial y^2} \leq 0 \end{cases} \tag{2.13}$$

We say that $f \in \mathcal{M}^{++}$ if $f \in \mathcal{M}^+$ and:

$$\begin{cases} \frac{\partial^2 f}{\partial x^2}, \frac{\partial^2 f}{\partial y^2} \leq 0 \\ \frac{\partial^3 f}{\partial x^2 \partial y}, \frac{\partial^3 f}{\partial x \partial y^2} \leq 0 \\ \frac{\partial^4 f}{\partial x^2 \partial y^2} \geq 0 \end{cases} \tag{2.13}$$

□

Using operator $K(s,t) = -(F(s,t) - F^X(s) - F^Y(t))$ we can define:

$$H(x,y,F) = \int_0^x \int_0^y F(s,t)ds dt \tag{2.14}$$

$$L(x,y,F) = \int_0^x \int_0^y K(s,t)ds dt \tag{2.15}$$

With upper indexes we denote its marginal versions:

$$H^X(x,F) = \int_0^x F^X(s)ds \tag{2.16}$$

$$H^Y(y,F) = \int_0^y F^Y(t)dt \tag{2.17}$$

and with lower indexes the operators for $F_1$ and $F_2$.

**Theorem 2.4.** (Second Order Bivariate Stochastic Dominance, Absolutely Continuous Case). Let $(X_1, Y_1)$ and $(X_2, Y_2)$ be random vectors with cdfs $F_1$ and $F_2$ respectively. Let's suppose in addition that both distributions are absolutely continuous with respect to Lebesgue measure and have common support with compact closure $\bar{U} = [0,1] \times [0,1]$.

1) Assume that:

$$\begin{cases} H_1^X(x) \leq H_2^X(x) & \forall x \in [0,1] \\ H_1^Y(y) \leq H_2^Y(y) & \forall y \in [0,1] \end{cases} \tag{2.18}$$

and that $H_1(x,y) \leq H_2(x,y)$ in $U$. Then $E(\phi(X_1,Y_1)) \geq E(\phi(X_2,Y_2))$ holds for every $C^4$ function $\phi \in \mathcal{M}^{--}$ which is both concave and increasing in each of its arguments.

2) Assume that:

$$\begin{cases} H_1^X(x) \leq H_2^X(x) & \forall x \in [0,1] \\ H_1^Y(y) \leq H_2^Y(y) & \forall y \in [0,1] \end{cases} \tag{2.19}$$

and that $L_1(x,y) \leq L_2(x,y)$ in $U$. Then $E(\phi(X_1,Y_1)) \geq E(\phi(X_2,Y_2))$ holds for every $C^4$ function $\phi \in \mathcal{M}^{++}$ which is both concave and increasing in each of its arguments.

**Proof.** Going back to (2.7) we have $E(\phi(X_1,Y_1)) = A_1 + B_1 + C_1 + \phi(1,1)$. Integration by parts and definition of $H_1^Y(y)$ gives:

$$B_1 = \frac{\partial \phi}{\partial y}(1,1) \int_0^1 F_1^Y(t)dt + \int_0^1 \frac{\partial^2 \phi}{\partial y^2}(1,t) H_1^Y(t)dt \tag{2.20}$$

Without loss of generality[14] we can assume $\phi$ has compact support in $U$ which cancels the first term. The assumed inequalities imply $B_1 \geq B_2$. Calculations to show $C_1 \geq C_2$ are obviously analogous.

These inequalities are valid for both (1) and (2) cases as we didn't use derivatives of order higher than 2. Terms $A_1$ and $A_2$ will separate cases.

Assume $\phi \in \mathcal{M}^{--}$ and put $\tilde{\phi}(x,y) = \frac{\partial^2 \phi}{\partial x \partial y}(x,y)$. By definition of $H_1(x,y)$ we have:

$$\int_U \frac{\partial^2 \phi}{\partial x \partial y}(s,t) F_1(s,t) ds dt = \int_U \tilde{\phi}(s,t) \frac{\partial^2 H_1}{\partial x \partial y}(s,t) ds dt \qquad (2.21)$$

This integral has the exact form of (2.1) with $\tilde{\phi}$ and $H_1$ replacing $\phi$ and $F_1$, so clearly:

$$\int_U \tilde{\phi}(s,t) \frac{\partial^2 H_1}{\partial x \partial y}(s,t) ds dt = \int_U H_1(s,t) \frac{\partial^2 \tilde{\phi}}{\partial x \partial y}(s,t) ds dt - \int_0^1 \frac{\partial \tilde{\phi}}{\partial y}(1,t) H_1^Y(t) dt +$$
$$\tilde{\phi}(1,1) - \int_0^1 \frac{\partial \tilde{\phi}}{\partial x}(s,0) H_1^X(s) ds \qquad (2.22)$$

Now we observe that $\frac{\partial \tilde{\phi}}{\partial x} = \frac{\partial^3 \phi}{\partial x^2 \partial y}$ and $\frac{\partial \tilde{\phi}}{\partial y} = \frac{\partial^3 \phi}{\partial y^2 \partial x}$ are both positive, and then this case follows with the same argument as in part (1) of Theorem 2.2.

Proof of part (2) of Theorem 2.4 starts with (2.12) and uses analogous properties for $\phi \in \mathcal{M}^{+\mp}$. □

A key point of both part (2) of Theorem 2.2 and Theorem 2.4 is that when we are dealing we supermodular classes we cannot ignore the effects of the boundary part of the integrals.

These effects show that in the supermodular case, the effect of joint distribution $F(x,y)$ on the expectation value can be compensated by the effects of marginal distributions and vice-versa. This turns out to be relevant in the interpretation of stochastic dominance within the framework of welfare economics, as shown by Perez[15].

---

[14] We can use density arguments in $C^4(U)$ if needed.

[15] Perez, L. (2015). Tests de Dominancia Estocástica en base a Estadísticos de Kolmogorov-Smirnov Multivariados, con Aplicaciones al Estudio de la Desigualdad Económica Multidimensional. PhD. Dissertation, UBA.

## 3. General Result on Modularity and Stochastic Dominance

Theorems 2.2 and 2.4, that is, both Hadar & Russell's (1974), Levy & Paroush's (1974) and Atkinson and Bourguignon's (1982) papers, assumed absolutely continuous distributions. This can be restrictive in applications as distributions are often assumed discrete, e.g. in the definition of inequality indexes in welfare economics.

We'll next show that both results continue to hold if we only assume that distributions have compact support. We may make use of the general Lebesgue-Stieljes integration by parts or even work in Sobolev spaces. But those techniques, although suitable for the general proof, would obscure the key point of our problem, the interaction between border and correlation effects in the supermodular case.

Our approach will be *direct* in the sense that we'll study Riemann-Stieljes sums and take limit. As Lebesgue-Stieljes integral can be seen as a completion of Riemann-Stieljes integral in the operator sense[16], the general result then follows.

**Theorem 3.1.** (First Order Bivariate Stochastic Dominance, General Case). The result in Theorem 2.2 holds for every pair of random vectors $(X_1, Y_1)$ and $(X_2, Y_2)$ cdfs $F_1$ and $F_2$ as long as they have compact support.

**Proof.** Without loss of generality we can assume $U = [0,1] \times [0,1]$.

For each distribution, we have to compute:

$$E(\phi(X,Y)) = \int_U \phi(s,t) dF(s,t) \qquad (3.1)$$

and compare the results. As functions $\phi$ involved are assumed to be increased and bounded on $U$, this expectations exist and are finite.

We take a partition $\Pi_n$:

$$\Pi_n = \{0 = x_0 < x_1 < \cdots < x_n = 1\} \times \{0 = y_0 < y_1 < \cdots < y = 1\} \qquad (3.2)$$

Inside each elementary block $B_{ij} = [x_{i-1}, x_i] \times [y_{j-1}, y_j]$ we select a point $(x_i^0, y_j^0)$ as in Fig.3.1. We are going to evaluate $\phi$ at these points.

Each elementary block $B_{ij}$ has a *quasi-volume* $\sigma(B_{ij})$ given by the double variation of the distribution function:

$$\sigma(B_{ij}) = F(x_i, y_j) + F(x_{i-1}, y_{j-1}) - F(x_i, y_{j-1}) - F(x_{i-1}, y_j) \qquad (3.3)$$

---

[16] Hewitt, E. (1960). Integration by Parts for Stieljes Integrals, Am. Math. Monthly, 67 (5), May, 419-423.

We construct the Riemann-Stieljes sum for our partition $\Pi_n$ and function $\phi$ evaluated at the selected points:

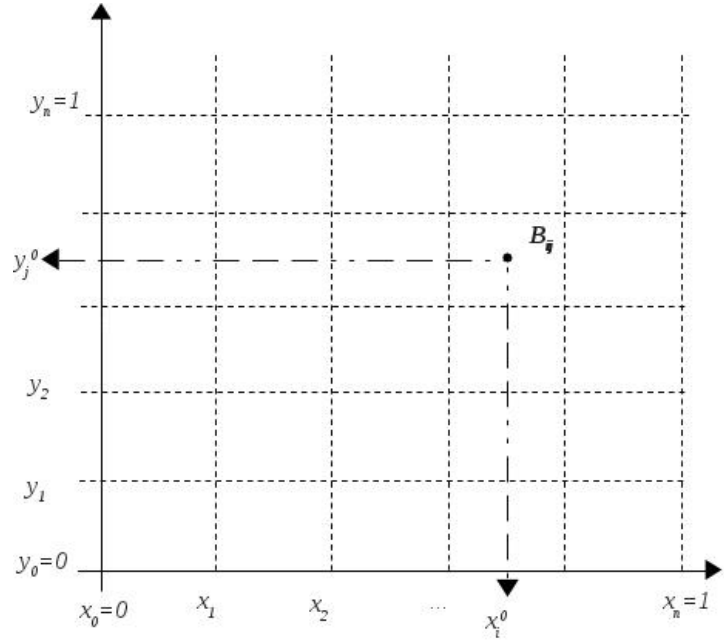

**Fig.3.1** Partition

$$S_{\Pi_n}(\phi) = \sum_{i,j=1}^{n} \phi(x_i^0, y_j^0)\sigma(B_{ij}) \qquad (3.4)$$

Or more explicitly:

$$S_{\Pi_n}(\phi) = \sum_{i,j=1}^{n} \phi\left(x_i^0, y_j^0\right)\left(F\left(x_i, y_j\right) + F\left(x_{i-1}, y_{j-1}\right) - F\left(x_i, y_{j-1}\right) - F\left(x_{i-1}, y_j\right)\right) \qquad (3.5)$$

For a better understanding of how the interior and border integrals impact on the expectation, we observe that if $(x_i, y_j)$ is not a border point of the partition (i.e. $x_i \neq 0; 1$ and $y_j \neq 0; 1$) then it appears at four terms in the sum, twice with positive and twice with negative sign (Fig.3.2).

The only points contributing to just two term in the sum are border points, i.e. points of the form $\{(x_i, 0), (x_i, 1), (0, y_j), (1, y_j)\}$. But as $F(s, t)$ is a distribution function with support $U$, we have $F(x_i, 0) = F(0, y_j) = 0$.

So the only relevant border effects are those of $F(x_i, 1)$ and $F(1, y_j)$ each appearing twice in the sum, once with positive and once with negative sign (Fig.3.3).

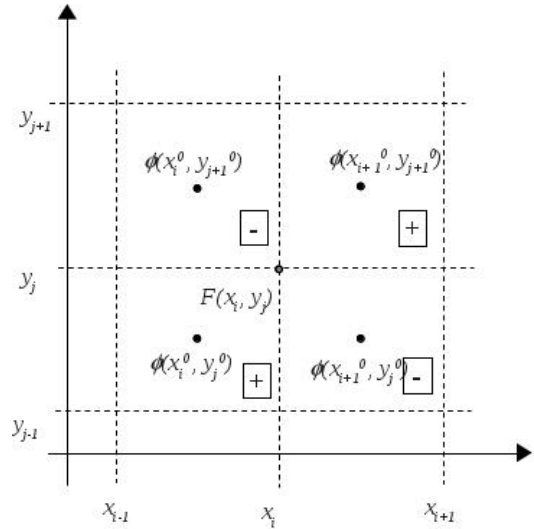

**Fig.3.2.** Contribution of interior terms to Stieljes sum.

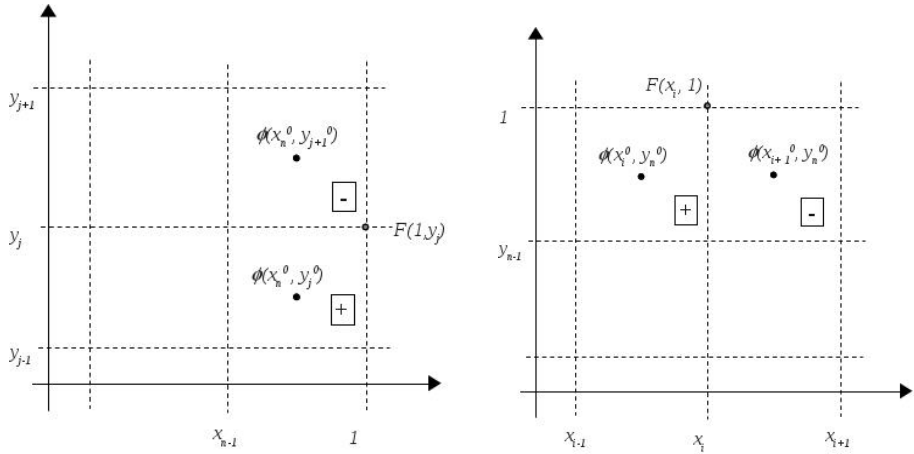

**Fig.3.3.** Effect of border terms in Stieljes sums.

We can gather terms including factor $F(x_i, y_j)$ for interior points:

$$F(x_i, y_j) * \left( \phi(x_i^0, y_j^0) + \phi(x_{i-1}^0, y_{j-1}^0) - \phi(x_{i-1}^0, y_j^0) - \phi(x_i^0, y_{j-1}^0) \right) =$$
$$F(x_i, y_j) \Delta_{ij} \phi \tag{3.6}$$

And for border points we have:

$$F(x_i, 1) * \left(\phi(x_i^0, y_n^0) - \phi(x_{i+1}^0, y_n^0)\right) = F(x_i, 1)\delta_{i,n}\phi \tag{3.7}$$

$$F(1, y_j) * \left(\phi(x_n^0, y_j^0) - \phi(x_n^0, y_{j+1}^0)\right) = F(1, y_j)\delta_{n,j}\phi \tag{3.8}$$

Finally, factor $F(1,1) = 1$ appears in only one term, so we can express:

$$S_{\Pi_n}(\phi) = \sum_{i,j=1}^{n-1} F(x_i, y_j)\Delta_{ij}\phi + \sum_{i=1}^{n-1} F(x_i, 1)\delta_{i,n}\phi + \sum_{j=1}^{n-1} F(1, y_j)\delta_{n,j}\phi + \phi(1,1) = A + B + C + \phi(1,1) \tag{3.9}$$

Note that in both the super and submodular case, we are assuming increasing $\phi$ so factors $\delta_{i,n}\phi$ and $\delta_{n,j}\phi$ are both negative.

Once we have computed Stieljes sums for each distribution we have:

$$S_{\Pi_n}^1 = A_1 + B_1 + C_1 + \phi(1,1) \tag{3.10}$$

$$S_{\Pi_n}^2 = A_2 + B_2 + C_2 + \phi(1,1) \tag{3.11}$$

Let's assume $\phi \in \mathcal{M}^-$ so $\frac{\partial^2 \phi}{\partial x \partial y}(x,y) \leq 0$ and then using Lagrange mean value theorem, it is easy to prove[17] that $\Delta_{ij}\phi \leq 0$. So if $F_1(s,t) \leq F_2(s,t)$ then obviously $A_1 \geq A_2$.

Taking limits it is clear that $F_1(x_i, 1) \leq F_2(x_i, 1)$ (same for the other limit) and as $\delta_{i,n}\phi$ and $\delta_{n,j}\phi$ are both negative, it is evident that $B_1 \geq B_2$ and $C_1 \geq C_2$. So we have, for every partition, $S_{\Pi_n}^1 \geq S_{\Pi_n}^2$ and this implies:

$$\int_U \phi(s,t) dF_1(s,t) \geq \int_U \phi(s,t) dF_2(s,t) \tag{3.12}$$

Or more succinctly:

$$E(\phi(X_1, Y_1)) \geq E(\phi(X_2, Y_2)) \tag{3.13}$$

for every $\phi \in \mathcal{M}^-$. This proves part (1).

Note that in this case (submodular) the effects of the interior and border parts of the sums (and integrals) are in the same direction, so there is no compensation between them.

Let's use the shorter notation $\phi(x_i^0, y_j^0) = \phi_{ij}$. For the supermodular case we can observe that:

---

[17] Topkis, D. (1998). *Supermodularity and Complementarity*. Princeton University Press.

$$\sum_{i=1}^{n-1} \Delta_{i,j}\phi = \phi_{1,j-1} - \phi_{1,j} + \phi_{n,j} - \phi_{n,j-1} = \phi_{1,j-1} - \phi_{1,j} - \delta_{n,j}\phi \Rightarrow$$

$$\Rightarrow \delta_{n,j}\phi = \phi_{1,j-1} - \phi_{1,j} - \sum_{i=1}^{n-1}\Delta_{i,j}\phi = \delta_{1,j}\phi - \sum_{i=1}^{n-1}\Delta_{i,j}\phi \qquad (3.14)$$

where we used the overlapping nature of the $\Delta_{i,j}\phi$ sequence to cancel intermediate terms, as shown in Fig.3.4.

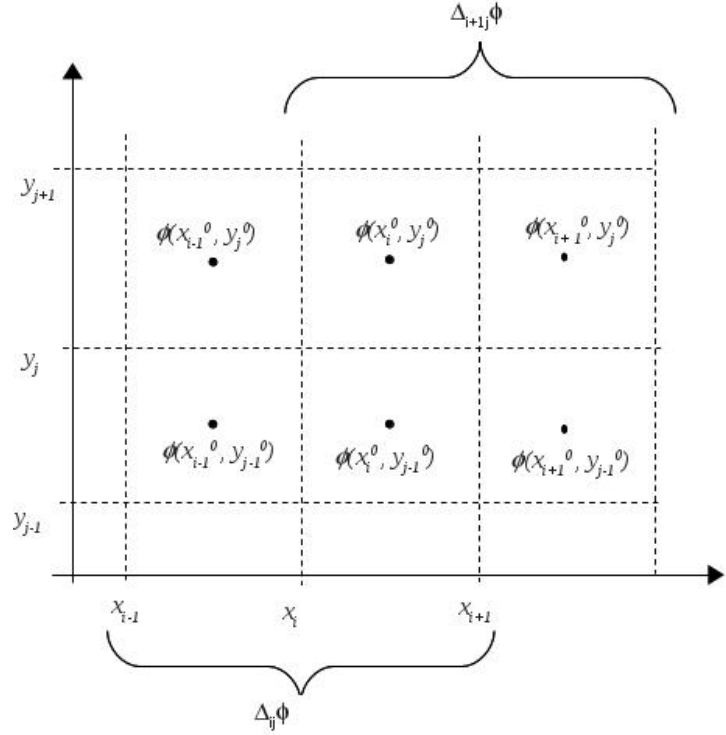

**Fig.3.4** Overlaping terms

Analogous calculations give us:

$$\delta_{i,n}\phi = \delta_{i,1}\phi - \sum_{j=1}^{n-1}\Delta_{i,j}\phi \qquad (3.15)$$

Replacing these in $B$ and $C$ we have:

$$B = \sum_{i=1}^{n-1} F(x_i,1)\delta_{i,1}\phi - \sum_{i,j=1}^{n-1} F(x_i,1)\Delta_{i,j}\phi \qquad (3.16)$$

$$C = \sum_{j=1}^{n-1} F(1,y_j)\delta_{1,j}\phi - \sum_{i,j=1}^{n-1} F(1,y_j)\Delta_{i,j}\phi \qquad (3.17)$$

Returning to (3.9) with all these results, we obtain:

$$S_{\Pi_n}(\phi) = \sum_{i,j=1}^{n-1}\left(F(x_i,y_j) - F(x_i,1) - F(1,y_j)\right)\Delta_{ij}\phi + \sum_{i=1}^{n-1} F(x_i,1)\delta_{i,1}\phi +$$
$$\sum_{j=1}^{n-1} F(1,y_j)\delta_{1,j}\phi + \phi(1,1) = -\sum_{i,j=1}^{n-1} K(x_i,y_j)\Delta_{ij}\phi + \sum_{i=1}^{n-1} F(x_i,1)\delta_{i,1}\phi +$$
$$\sum_{j=1}^{n-1} F(1,y_j)\delta_{1,j}\phi + \phi(1,1) \qquad (3.18)$$

Both $\delta_{1,j}\phi$ and $\delta_{i,1}\phi$ are negative as $\phi$ is increasing in each variable and for the supermodular case we have $\Delta_{ij}\phi \geq 0$.

Taking this into account and recalling assumptions $F_1^X(s) \leq F_2^X(s)$, $F_1^Y(t) \leq F_2^Y(t)$ and $K_1(s,t) \leq K_2(s,t)\ \forall (s,t) \in U$, we have as before:

$S_{\Pi_n}^1 \geq S_{\Pi_n}^2$ and we know this implies:

$$E(\phi(X_1,Y_1)) \geq E(\phi(X_2,Y_2)) \qquad (3.19)$$

And this completes the proof. $\square$

Note that this proof, in addition to extending Hadar & Russell's and Levy & Paroush's results to the general (not necessarily continuous) case gives us a *formula* for the discrete case, obtained by suitably reinterpreting Stieljes sums.

As in the case for first order dominance, we observe that an inequality like $E(\phi(X_1,Y_1)) \geq E(\phi(X_2,Y_2))$ and the sufficient conditions found in Theorems 2.2 and 2.4, make no explicit use of the absolute continuity of the underlying distributions. This continuity property only arose in the *proof* of the result, as it provided an easy way to go using integration by parts.

It has to be clear, though, that the conditions for stochastic dominance hold for general compact support distributions.

We turn now our attention to second order dominance.

**Theorem 3.2.** (Second Order Bivariate Stochastic Dominance, General Case). The result in Theorem 2.4 holds for every pair of random vectors $(X_1,Y_1)$ and $(X_2,Y_2)$ cdfs $F_1$ and $F_2$ as long as they have compact support.

**Proof.** We are not going to give a complete proof, but only to show that appropriate reinterpretation of terms renders us the same conditions as in the proof of Theorem 3.1.

Going back to (3.9) we have $A = \sum_{i,j=1}^{n-1} F(x_i,y_j)\Delta_{ij}\phi$. Now, observing that:

$$\Delta_{ij}\phi = \phi(x_i^0, y_j^0) + \phi(x_{i-1}^0, y_{j-1}^0) - \phi(x_{i-1}^0, y_j^0) - \phi(x_i^0, y_{j-1}^0) \qquad (3.20)$$

and using Lagrange mean value theorem, we have:

$$\Delta_{ij}\phi = \frac{\partial^2 \phi}{\partial x \partial y}(\tilde{x}_i^0, \tilde{y}_j^0)\Delta x \Delta y + \frac{\partial^2 \phi}{\partial x^2}(\hat{x}_i^0, \hat{y}_j^0) O(\Delta x)^2 \qquad (3.21)$$

where $(\tilde{x}_i^0, \tilde{y}_j^0)$ and $(\hat{x}_i^0, \hat{y}_j^0)$ are suitable points in the rectangle $[x_{i-1}^0, x_i^0] \times [y_{j-1}^0, y_j^0]$ and where we called $\Delta x = x_i^0 - x_{i-1}^0$ and $\Delta y = y_j^0 - y_{j-1}^0$.

Introducing this factor in $A$ yields:

$$\sum_{i,j=1}^{n-1} F(x_i, y_j) \Delta_{ij}\phi = \sum_{i,j=1}^{n-1} F(x_i, y_j) \frac{\partial^2 \phi}{\partial x \partial y}(\tilde{x}_i^0, \tilde{y}_j^0) \Delta x \Delta y +$$
$$\sum_{i,j=1}^{n-1} F(x_i, y_j) \frac{\partial^2 \phi}{\partial x^2}(\hat{x}_i^0, \hat{y}_j^0) O(\Delta x)^2 \quad (3.22)$$

Now the quadratic factor in the last sum implies this term will converge to 0 as the diameter of the partition tends to 0, so this term is not of interest.

Now for the first ter we can observe that $F(x_i, y_j)\Delta x \Delta y \approx \Delta_{ij} H$ where exact equality holds if instead of $H$ we take the Riemann sums of the double integral that defines it.

So we know that the limit of sums in $\sum_{i,j=1}^{n-1} F(x_i, y_j) \frac{\partial^2 \phi}{\partial x \partial y}(\tilde{x}_i^0, \tilde{y}_j^0) \Delta x \Delta y$ is the same as that of the sums in $\sum_{i,j=1}^{n-1} \tilde{\phi}(x_i, y_j) \Delta_{ij} H$. But this is just the same as starting over the proof of Theorem 3.1. with $\tilde{\phi}$ and $H$ instead of $\phi$ and F.

It should be clear from this point that using the inequalities assumed for $H$ and for the marginal distributions will give us the desired result. □

## 4. Discussion and Concluding Remarks

The results we obtained in Section 3 are relevant in the empirical implementation of the concept of stochastic dominance, as we can use the sufficient conditions for dominance without assuming absolute continuity of the underlying distributions.

The importance of this generalization of the previous results lies on the fact that absolute continuity is too restrictive in two aspects that are essential in most applications.

First, assuming absolute continuity of the underlying distributions leaves out the case of discrete or even finite distributions. As finite and discrete distributions are very common in the empirical literature on welfare economics and economic inequality measurement[18], extending sufficient conditions for multivariate stochastic dominance from the absolutely continuous case to the general compact support case is a significant contribution to the field of applied quantitative economics.

---

[18] Cowell, S. (2001). Estimation of Inequality Indices. In Silber, J. (Ed.), Handbook of Income Inequality Measurement, New York: Springer Science+Business Media.

On the second hand, for statistical implementations, the extension to the general compact support case is relevant when taken together with Kolmogorov-Smirnov-type test statistics. Although the usual implementations of Kolmogorov-Smirnov statistics assume continuity of the underlying distribution, in the most general case we can use empirical processes theory and bootstrap techniques to obtain limit distributions, as done by van der Waart & Wellner (2000). Using this theoretical results, the work of Perez (2015) finds weak limits for test statistics based on the sufficient conditions proven in the present paper.

Regarding the complexity of the proof, both in the continuous and the general case, we observe that it arises from the *correlation* of the dimensions. The paper by Scarsini (1988)[19] provides sufficient conditions for bivariate stochastic dominance in the case were *X* and *Y* are independent. The resulting proof is much easier but the relevance of the result is strongly diminished by the fact that in almost any conceivable economic or financial application the assumption of no correlation is unacceptable. Atkinson & Bourguignon (1982), for example, study bivariate economic inequality using the dimensions of income and health and Perez (2015) uses income and human capital. In both cases, there are both theoretical and empirical reasons to rule out independence of the variables.

Another trick to find easier proofs for bivariate dominance conditions is to consider only bivariate distributions with identical marginal distributions, as done by Tchen (1980)[20]. Evidently, as we saw that in the supermodular cases the contribution of the marginal distributions add to that of the bivariate cdf (or its integrals), assuming identical marginals will ease the proof and give simpler dominance conditions. But, then again, there is no theoretical or empirical reasons to assume that this is the case in the usual applications. It has to be clear that aside for the inequalities assumed in the Theorems, our proof holds for arbitrary marginals, making then for a meaningful extension of the literature.

The result in Epstein & Tanny (1980)[21] using *increasing correlation transformations* has a meaningful interpretation in welfare economics, as it acts as analogue to Pigou transfer principle. But as the concept has no direct statistical implementation, at least not in terms of easy transformations of underlying cdfs, its applications are mostly theoretical. In this way, a *direct*

---

[19] Scarsini, M. (1988). Multivariate Stochastic Dominance With Fixed Dependence Structure, *Technical Report N.255*, Department of Statistics, Stanford Univesity.

[20] Tchen, A. (1980). Inequalities for Distributions with Given Marginals, *The Annals of Probability, 8 (4)*, 814-827.

[21] Epstein, L. & Tanny, S. (1980). Increasing Generalized Correlation: A Definition and Some Economic Consequences, *The Canadian Journal of Economics / Revue Canadienne d'Economique, 13 (1)*, February, 16-34.

*approach* to stochastic dominance as taken in the present paper is more useful in applications.

We emphasize the importance of the concept of multivariate stochastic dominance in the welfare economics literature as it gives a way to compare multivariate prospects with direct economic interpretation as it involves expected values of utility functions. In this regard, the mentioned *direct approach* to stochastic dominance that we use is of fundamental importance. *Inverse* stochastic dominance, including Lorenz dominance, although useful in the univariate case, has no meaning in the multivariate case as there is no univocal direct analogue to Lorenz curve in *n* dimensions[22]. So if we have to deal with multidimensional economic inequality, we have to use the direct approach given in the present paper.

There are two natural extensions for the results proven in this paper. First, we can extend the results to higher dimensions. Second, we can prove analogous results for higher orders of dominance. Formally, there is no difficulty in doing this, although two points of warning have to be made, concerning to the limitations of our approach.

Regarding the first point, even for dimension n=2, as we study higher order dominances, the inequalities involved can present complications. As we generalize the integral operators of Davidson & Duclos (2000) that were given in equations (1.3) and (1.4), we see that for each order of dominance we have to deal with a new double integral. Working as in our proof, or even for the absolutely continuous case, we see that each use of integration by parts formula will introduce higher derivatives of the test functions $\phi$.

Thus the structure of modularity classes will become more and more complex as more conditions on the sign of those derivatives are needed. This has the problem that some of the conditions involved in the definition of modularity classes may not have direct interpretation in economics of risk analysis, thus implying some arbitrariness in the selected classes. This is the main reason to restrict the analysis of stochastic dominance to first and second order.

The second point is perhaps more relevant in applications, as it involves even first and second order dominance. As the number of points in data sets needed to obtain a given level of statistical significance grows exponentially with the dimension of the distributions[23], the extension of our results to the higher dimensional case will undoubtedly fall in the *curse of dimensionality.*

---

[22] Perez, L. (2015). Tests de Dominancia Estocástica en base a Estadísticos de Kolmogorov-Smirnov Multivariados, con Aplicaciones al Estudio de la Desigualdad Económica Multidimensional. PhD. Dissertation, UBA.

[23] Devroye, L. & Lugosi, G. (2001). Combinatorial Methods in Density Estimation. New York: Springer.

Computational methods as the bootstrap will then be costly to implement, restricting the applications of our results.

Despite these restrictions, our results are relevant for the applied literature as they cover the bivariate case where involved dimensions are income and human capital, making for meaningful welfare interpretation of the dominance results.

Interesting contributions to stochastic dominance analysis and its applications in finance and economics continue to appear to this day, as shown e.g. by Meyer & Strulovici (2012) or the complete Issue 4 of Volume 147 of the Journal of Economic Theory, appropriately labelled *Inequality and Risk*. The present paper is intended as a contribution to this growing literature.

**References.**